\documentclass[preprint,11pt]{elsarticle}
\usepackage{amsmath,amssymb,latexsym}

\newtheorem{theorem}{Theorem}
\newtheorem{lemma}{Lemma}[section]

\newtheorem{corollary}[lemma]{Corollary}
\newtheorem{proposition}[lemma]{Proposition}

\newtheorem{definition}[lemma]{Definition}
\newtheorem{example}[lemma]{Example}
\newtheorem{remark}[lemma]{Remark}

\newcommand{\bl}{\begin{lemma}}
\newcommand{\el}{\end{lemma}}
\newcommand{\bt}{\begin{theorem}}
\newcommand{\et}{\end{theorem}}
\newcommand{\bcor}{\begin{corollary}}
\newcommand{\ecor}{\end{corollary}}
\newcommand{\bp}{\proof{.}}
\newcommand{\ep}{\eop}
\newcommand{\bpr}{\begin{proposition}}
\newcommand{\epr}{\end{proposition}}
\newcommand{\brem}{\begin{remark} \em}
\newcommand{\erem}{\end{remark}}
\newcommand{\bd}{\begin{definition} \em}
\newcommand{\ed}{\end{definition}}
\newcommand{\bex}{\begin{example} \em
}
\newcommand{\eex}{\end{example}}
\newcommand{\beq}{\begin{equation} }
\newcommand{\eeq}{\end{equation}}

\newcommand{\bi}{\begin{itemize}
  }
\newcommand{\ei}{\end{itemize}}
\newcommand{\ben}{\begin{enumerate} }
\newcommand{\een}{\end{enumerate} }

\newcommand{\refeq}[1]{(\ref{#1})}
\newenvironment{enumr}{
\renewcommand{\theenumi}{\roman{enumi}}
\renewcommand{\labelenumi}{$\mathrm{(\theenumi)}$}
\begin{enumerate}     }{\end{enumerate}
\renewcommand{\theenumi}{\arabic{enumi}}
\renewcommand{\labelenumi}{\theenumi.}}
\newcommand{\benr}{\begin{enumr}
  }
\newcommand{\eenr}{
\end{enumr}}

\newcommand{\bs}{\bigskip}

\newcommand{\ul}[1]{\underline{#1}}
\newcommand{\al}[1]{\forall #1\:}
\newcommand{\ex}[1]{\exists #1\:}

\newlength{\hilflh}

\newcommand{\naturals}{\mathbb{N}}

\renewcommand{\emptyset}{\varnothing}

\newcommand{\cL}{{\mathcal L}}

\newcommand{\cV}{{\mathcal V}}

\newcommand{\cW}{{\mathcal W}}

\newcommand{\cC}{{\mathcal C}}

\newcommand{\ga}{\alpha}
\newcommand{\gb}{\beta}

\renewcommand{\ge}{\varepsilon}

\newcommand{\gs}{\sigma}
\newcommand{\gy}{\gamma}
\newcommand{\gw}{\omega}
\newcommand{\gS}{\Sigma}

\renewcommand{\phi}{\varphi}

\newcommand{\eqv}{\leftrightarrow}

\renewcommand{\Pr}{\mathrm{Pr}}

\newcommand{\ol}{\overline}

\newcommand{\GL}{\mathbf{GL}}
\newcommand{\GLP}{\mathbf{GLP}}

\newcommand{\Con}{\mathrm{Con}}

\newcommand{\Prf}{\mathrm{Prf}}

\newcommand{\RFN}{\mathsf{RFN}}

\newcommand{\PA}{\mathsf{PA}}

\newcommand{\rst}{\upharpoonright}
\newcommand{\gn}[1]{\ulcorner #1 \urcorner}

\newcommand{\fc}{\Vdash}      %amsforcing
\renewcommand{\models}{\vDash}      %amsmodels
\newcommand{\nfc}{\nVdash}

\newcommand{\Imp}{\Rightarrow}

\newcommand{\Var}{\mathrm{Var}}

\newcommand{\nat}{\naturals}

\newcommand{\J}{\mathbf{J}}

\renewcommand{\leq}{\leqslant}
\renewcommand{\geq}{\geqslant}

\newcommand{\RCo}{\mathbf{RC\omega}}
\newcommand{\Rc}{\mathrm{RC}}
\newcommand{\rj}{\mathrm{RJ}}

\newcommand{\eop}{$\Box$ \protect\par \addvspace{\topsep}}
\newcommand{\proof}[1]{\protect\par\addvspace{\topsep}\noindent {\bf Proof#1}}
\newcommand{\RC}{\mathbf{RC}}
\newcommand{\RJ}{\mathbf{RJ}}

\journal{Annals of Pure and Applied Logic}

\begin{document}

\begin{frontmatter}
\title{Positive provability logic \\ for uniform reflection principles}

\author{Lev Beklemishev}
\ead{bekl@mi.ras.ru} \ead[url]{http://www.mi.ras.ru/~bekl}

\address{Steklov Mathematical Institute\fnref{aff},
Russian Academy of Sciences,\\ Gubkina str. 8, 119991, Moscow,
Russia} \fntext[aff]{Also affiliated at: Moscow M.V. Lomonosov State
University and National Research University Higher School of
Economics.}

\begin{abstract}
We deal with the fragment of modal logic consisting of implications
of formulas built up from the variables and the constant `true' by
conjunction and diamonds only. The weaker language allows one to
interpret the diamonds as the uniform reflection schemata in
arithmetic, possibly of unrestricted logical complexity. We
formulate an arithmetically complete calculus with modalities
labeled by natural numbers and $\gw$, where $\gw$ corresponds to the
full uniform reflection schema, whereas $n<\gw$ corresponds to its
restriction to arithmetical $\Pi_{n+1}$-formulas. This calculus is
shown to be complete w.r.t.\ a suitable class of finite Kripke
models and to be decidable in polynomial time.
\end{abstract}


\end{frontmatter}

\section{Introduction}
Several applications of provability logic in proof theory made use
of a polymodal logic $\GLP$ due to Giorgi
Japaridze~\cite{Dzh86,Boo93}. This system, although decidable, is
not very easy to handle. In particular, it is not Kripke
complete~\cite{Boo93}. It is complete w.r.t.\ the more general
topological semantics, however this could only be established
recently by rather complicated techniques \cite{BekGab11}.

A weaker system, called \emph{Reflection Calculus} and denoted
$\RC$, was introduced in \cite{Bek12a}. It is much simpler than
$\GLP$ yet expressive enough to regain its main proof-theoretic
applications. It has been outlined in \cite{Bek12a} that $\RC$
allows to define a natural system of ordinal notations up to $\ge_0$
and serves as a convenient basis for a proof-theoretic analysis of
Peano Arithmetic in the style of \cite{Bek04,Bek05}. This includes a
consistency proof for $\PA$ based on transfinite induction up to
$\ge_0$, a characterization of its $\Pi_n^0$-consequences in terms
of iterated reflection principles, and a combinatorial independence
result.

From the point of view of modal logic, $\RC$ can be seen as a
fragment of polymodal logic consisting of implications of the form
$A\to B$, where $A$ and $B$ are formulas built-up from $\top$ and
propositional variables using just $\land$ and the diamond
modalities. We call such formulas $A$ and $B$ \emph{strictly
positive} and will often omit the word
`strictly.'\footnote{Traditionally, positive modal formulas may also
involve disjunctions and box modalities. However, in the present
paper we will not consider positive formulas in this more general
sense.}

A somewhat different but equivalent axiomatization of $\RC$ (as an
equational calculus) has been earlier found by Evgeny Dashkov in his
paper~\cite{Das12} which initiated the study of strictly positive
fragments of provability logics. Dashkov proved two important
further facts about $\RC$ which sharply contrast with the
corresponding properties of $\GLP$. Firstly, $\RC$ is complete with
respect to a natural class of finite Kripke frames. Secondly, $\RC$
is decidable in polynomial time, whereas most of the standard modal
logics (including $\GL$ and $\GLP$) are \textsc{PSpace}-complete.

Another advantage of going to a strictly positive language is
explored in the present paper. Strictly positive modal formulas
allow for more general arithmetical interpretations than those of
the standard modal logic language. In particular, propositional
formulas can now be interpreted as arithmetical \emph{theories}
rather than individual \emph{sentences}. (Notice that the `negation'
of a theory would not be well-defined.)

Similarly, the diamonds need no longer be interpreted as individual
\emph{consistency assertions} but as the more general
\emph{reflection schemata} not necessarily having finite
axiomatizations. Thus, for example, the full uniform reflection
schema can be considered as a modality in the context of positive
provability logic (see \cite{KrL,Bek05} for general information on
reflection principles). Such interpretations are not only natural
but can be useful for further development of the approach to
proof-theoretic analysis via provability algebras. Thus, positive
provability logic allows to speak about certain notions not nicely
representable in the context of the standard modal logic.

The main contribution of this paper is a Solovay-style arithmetical
completeness result for an extension of $\RC$ by a new modality
corresponding to the unrestricted uniform reflection principle. This
is the primary example of a modality not representable in the full
modal logic language. The system obtained is shown to be decidable
and to enjoy a suitable complete Kripke semantics along with the
finite model property.

Whereas the modal logic part of our theorem is a simple extension of
Dashkov's results, the arithmetical part is more substantial. We
introduce a new modification of the Solovay construction using some
previous ideas from \cite{Dzh86, Ign93, Bek11}. Since the
arithmetical complexity of the unform reflection schema is
unbounded, a single Solovay-style function is not enough for our
purpose. Instead, we deal with infinitely many Solovay functions, of
increasing arithmetical complexity, uniformly and simultaneously.%
\footnote{An interesting arithmetical completeness proof for an
extension of modal logic $\GLP$ by transfinitely many modalities has
recently appeared in \cite{FJ13}. However, the considered
interpretation is different and not applicable in our situation.}

The paper is organized as follows. Firstly, we introduce positive
modal language and the systems leading to the arithmetically
complete reflection calculus $\RCo$. Secondly, we present the
details of its arithmetical interpretation and somewhat tediously
prove the corresponding soundness theorem. Thirdly, we study the
Kripke semantics of positive provability logics and obtain
completeness results, along with a suitable version of the finite
model property. Fourthly, we obtain polynomial complexity bounds for
the derivability problem in $\RCo$ by adapting the techniques of
Dashkov. Finally, we prove the main result of this paper, the
arithmetical completeness theorem for $\RCo$.

\section{Reflection calculus and its basic properties}

Consider a modal language $\cL$ with propositional variables
$p,q$,\dots , a constant $\top$ and connectives $\land$ and $\ga$,
for each ordinal $\ga\leq\gw$ (understood as diamond modalities).
Strictly positive formulas (or simply \emph{formulas}) are built up
by the grammar:
$$A::= p \mid \top \mid (A\land B) \mid \ga A, \quad \text{where $\ga\leq\gw$.}$$
\emph{Sequents} are expressions of the form $A\vdash B$ where $A,B$
are strictly positive formulas. The system $\RJ$ is given by the
following axioms and rules:

\ben
\item $A\vdash A; \quad A\vdash\top; \quad$ if $A\vdash B$ and $B\vdash C$ then $A\vdash C$ (syllogism);
\item $A\land B\vdash A; \quad A\land B\vdash B; \quad$ if $A\vdash B$ and $A\vdash C$ then
$A\vdash B\land C$;
\item  if $A\vdash B$ then $\ga A\vdash \ga B$;\quad $\ga\ga A\vdash \ga A$;
\item $\ga\gb A\vdash \gb A;\quad$ $\gb\ga A\vdash \gb A$ for
$\ga\geq \gb$;
\item $\ga A\land \gb B\vdash \ga(A\land \gb B)$ for
$\ga >\gb$. \een

The systems $\RC$ and $\RCo$ are obtained from $\RJ$ by adding
respectively one or two of the following principles:

\ben\setcounter{enumi}{5}
\item $\ga A\vdash \gb A \text{ for $\ga>\gb$};$ \quad (monotonicity)
\item $\gw A \vdash A$. \quad (persistence)
\een

Dashkov~\cite{Das12} showed that $\RJ$, restricted to the language
without $\gw$ modality, axiomatizes the strictly positive fragment
of the polymodal logic $\J$ \cite{Bek10}, whereas $\RC$ axiomatizes
the strictly positive fragment of $\GLP$.

Notice that Axioms 4 are redundant in the presence of Axiom 6: if
$\ga \geq \gb$ then $\ga\gb A\vdash \gb\gb A\vdash \gb A$ and
$\gb\ga A\vdash \gb\gb A\vdash \gb A$.

If $L$ is a logic, we write $A\vdash_L B$ for the statement that the
sequent $A\vdash B$ is provable in $L$. As a simple example,
consider the sequent $$\gw (p\land q) \vdash (\gw p \land \gw q).$$
It is provable in $\RJ$ as follows: We have $p\land q \vdash p$,
hence $\gw(p\land q)\vdash \gw p$. Similarly, $\gw(p\land q)\vdash
\gw q$, therefore $\gw(p\land q)\vdash (\gw p\land \gw q)$, by the
conjunction rules. In contrast, $(\gw p\land \gw q) \nvdash_\RCo
\gw(p\land q)$, as we shall see below by a simple Kripke model
argument.

Formulas $A$ and $B$ are called \emph{$L$-equivalent} (written
$A\sim_L B$) if $A\vdash_L B$ and $B\vdash_L A$.

We also consider the fragments of various logics obtained by
restricting the language to a subset of modalities. Such a subset
$S\subseteq\gw+1=\{0,1,\dots,\gw\}$ is called a \emph{signature}. We
denote by $\cL_S$ the set of all strictly positive formulas in $S$.
Similarly, for a logic $L$ in $\cL$ we denote by $L_S$ the
restriction of the axioms and rules of $L$ to the language $\cL_S$.

For a positive formula $A$, let $\ell(A)$ denote $\{\ga\leq\gw:
\text{$\ga$ occurs in $A$}\}$. \bl \benr \item If $A\vdash_\RJ B$
then $\ell(B)\subseteq \ell(A)$;
\item If $A\vdash_L B$ where $L$ is $\RC$ or $\RC\gw$, then $\ell(B)\subseteq [0,\max\ell(A)]$.
\eenr \el

\bp\ In each case, this is proved by an easy induction on the length
of the derivation of $A\vdash B$. \ep

Let $C[A/p]$ denote the result of replacing in $C$ all occurrences
of a variable $p$ by $A$. If a logic $L$ contains Axioms 1, 2 and
the first part of 3, then $\vdash_L$ satisfies the following
\emph{positive replacement lemma}.

\bl Suppose $A\vdash_L B$, then $C[A/p]\vdash_L C[B/p]$, for any
$C$. \el \bp\ Induction on the build-up of $C$. \ep

A positive logic $L$ is called \emph{normal} if it contains the
rules 1, 2, and the first part of 3, and is closed under the
following \emph{substitution rule}: if $A\vdash_L B$ then
$A[C/p]\vdash_L B[C/p]$. It is clear that $\RJ$, $\RC$ and $\RC\gw$,
as well as their restricted versions, are normal.

\section{Arithmetical interpretation}

We define the intended arithmetical interpretation of the positive
modal language. The idea is that propositional variables (and
positive formulas) now denote possibly infinite theories rather than
individual sentences. To avoid possible problems with the
representation of theories in the language of $\PA$, we deal with
\emph{primitive recursive numerations} of theories rather than with
the theories as sets of sentences.

All theories in this paper will be formulated in the language of
Peano Arithmetic $\PA$ and contain the axioms of $\PA$. It is
convenient to assume that the language of $\PA$ contains the symbols
for all primitive recursive programs. A \emph{primitive recursive
numeration} of a theory $S$ is a bounded arithmetical formula
$\gs(x)$ defining the set of G\"odel numbers of the axioms of $S$ in
the standard model of arithmetic. Given such a $\gs$, we have a
standard arithmetical $\Sigma_1$-formula $\Box_\gs(x)$ expressing
the provability of $x$ in $S$ (see \cite{Fef60}). We often write
$\Box_\gs\phi$ for $\Box_\gs(\gn{\phi})$. The expression $\bar n$
denotes the numeral $0'\cdots'$ ($n$ times). If $\phi(v)$ contains a
parameter $v$, then $\Box_\gs\phi(\bar x)$ denotes a formula (with a
parameter $x$) expressing the provability of the sentence $\phi(\bar
x/v)$ in $S$.

Given two numerations $\gs$ and $\tau$, we write $\gs\vdash_{\PA}
\tau$ if $$\PA\vdash \al{x}(\Box_\tau(x)\to \Box_\gs(x)).$$  We
write $\gs\vdash \tau$ if $\nat\models \al{x}(\Box_\tau(x)\to
\Box_\gs(x)),$ that is, if the theory numerated by $\gs$ contains
the one numerated by $\tau$. We will only consider the numerations
$\gs$ such that $\gs\vdash_\PA \gs_\PA,$ where $\gs_\PA$ is some
standard numeration of $\PA$.

With any finite extension of $\PA$ of the form $\PA+\phi$ we will
associate its standard numeration $\gs_\PA\lor (x=\gn{\phi})$ that
will be denoted $\ul\phi$. For obvious reasons we have:
$\ul\phi\vdash_\PA \ul\psi$ iff $\ul\phi\vdash \ul\psi$ iff
$\PA+\phi\vdash\psi$. (The statement $\ul\phi\vdash_\PA \ul\psi$
implies $\PA+\phi\vdash \psi$ by the soundness of $\PA$, the
converse is formalizable in $\PA$.)

Given a numeration $\gs$ of $S$, the consistency of $S$ is expressed
by $\Con(\gs):=\neg\Box_\gs\bot$. A theory $S$ is called
\emph{$n$-consistent} if $S$ together with the set of all true
$\Sigma_{n+1}$-sentences is consistent. The $n$-consistency of $S$
is expressed by the formula
$$\al{x\in\Pi_{n+1}}(\Box_\gs (x)\to T_n(x)), \leqno \Con_n(\gs):$$
where $T_n$ is the standard $\Pi_{n+1}$-truthdefinition for
$\Pi_{n+1}$-formulas (see \cite{HP}) and $x\in\Pi_{n+1}$ denotes the
primitive recursive formula expressing that $x$ is a G\"odel number
of a $\Pi_{n+1}$-sentence.

Concerning the truthdefinitions we assume that $\PA\vdash \phi\eqv
T_n(\gn{\phi}),$ for each $\Pi_{n+1}$-sentence $\phi$. Moreover,
this very fact can be formalized in $\PA$ uniformly in $n$:
\beq\label{tr} \PA\vdash \al{n}\al{x\in\Pi_{n+1}}\Box_\PA(x\eqv
T_n(\bar x)),\eeq as the sequence of formulas $\gn{T_n}$ is
primitive recursive in $n$ and the corresponding proofs are
constructed inductively.

The formula $\Con_n(\gs)$ is often called \emph{the global
$\Pi_{n+1}$-reflection principle for $S$} and is denoted
$\RFN_{\Pi_{n+1}}(S)$ (see \cite{Smo77,Bek95}). We note that the
formula $\Con _0(\gs)$ is $\PA$-provably equivalent to $\Con(\gs)$.

The \emph{uniform reflection principle} for $S$ is the schema
$$\{\Con_n(\gs):n\in\gw\}. \leqno \Con_\gw(\gs):$$
It is well-known that $\Con_\gw(\gs)$ is $\PA$-provably equivalent
to the schema $$\al{x}(\Box_\gs \phi(\bar x)\to \phi(x)),$$ for each
arithmetical formula $\phi(x)$, which is usually denoted $\RFN(S)$.

The uniform reflection principle is elementarily axiomatized, and we
fix a standard function mapping any numeration $\gs$ to the
numeration of $\PA+\Con_\gw(\gs)$ (denoted $\ul{\Con}_\gw(\gs)$).
Similarly, the formula
$$\gs_\PA(x)\lor x=\gn{\Con_n(\gs)}$$ numerating the theory
$\PA+\Con_n(\gs)$ will be denoted $\ul \Con_n(\gs)$.

The intended arithmetical interpretation maps positive modal
formulas to primitive recursive numerations in such a way that
$\top$ corresponds to the standard numeration of $\PA$, $\land$
corresponds to the union of theories, $n$ corresponds to the
standard numeration of $\Con_n$, for each $n<\gw$, and $\gw$ to the
standard numeration of $\Con_\gw$.

\bd An \emph{arithmetical interpretation} is a map $*$ from positive
modal formulas to numerations satisfying the following conditions:

\bi
\item $\top^*=\gs_{\PA}$; \quad $(A\land B)^*=(A^*\lor B^*)$;
\item $(n A)^*=\ul{\Con}_n(A^*)$;
\quad $(\gw A)^*=\ul{\Con}_{\gw}(A^*)$. \ei \ed

It is clear that the value $A^*$ is completely determined by the
interpretations $p_1^*,\dots,p_n^*$ of all the variables occurring
in $A$.

\bpr[soundness] \label{sound} Suppose $A\vdash_\RCo B$, then
$A^*\vdash_{\PA} B^*$, for all arithmetical interpretations $*$.
\epr

\bp\ Induction on the length of proof of $A\vdash B$ in $\RCo$. The
validity of the first two groups of rules of $\RCo$ is obvious. We
treat the modal axioms and rules.

If $\gs\vdash_\PA \tau$ then clearly $\ul{\Con}_{n}(\gs)\vdash_{\PA}
\ul{\Con}_{n}(\tau)$, for each $n<\gw$. Since this fact is
formalizable in $\PA$, we also obtain
$\ul{\Con}_\gw(\gs)\vdash_{\PA} \ul{\Con}_\gw(\tau)$. Also, the
validity of the monotonicity axioms 6 is clear. Next we need the
following lemma.

\bl \label{reflc} \begin{enumr} \item Let $S$ be numerated by $\gs$
and $\phi\in\Pi_{n+1}$. If $S\vdash\phi$ then
$\PA+\Con_n(\gs)\vdash\phi$;
\item Statement (i) holds provably in $\PA$ uniformly in $n$, that is,
$$\PA\vdash \al{n}\al{x\in\Pi_{n+1}}(\Box_\gs(x)\to
\Box_{\ul{\Con}_n(\gs)}(x)).$$ \end{enumr} \el

\bp\ We only prove Statement (ii). We reason in $\PA$ as follows.

Assume $x\in\Pi_{n+1}$ and $\Box_\gs(x)$. Then $\Box_\PA(\bar
x\in\Pi_{n+1}\land \Box_\gs(\bar x))$. On the other hand, by the
definition of $\Con_n(\gs)$ $$\Box_{\ul \Con_n(\gs)}\al{y}(\Box_\gs
(y)\land y\in \Pi_{n+1}\to T_n(y)).$$ This yields $$\Box_{\ul
\Con_n(\gs)}(\Box_\gs (\bar x)\land \bar x\in \Pi_{n+1}\to T_n(\bar
x)),
$$ so we obtain $\Box_{\ul \Con_n(\gs)} T_n(\bar x)$, and hence $\Box_{\ul \Con_n(\gs)}
(x)$ by \refeq{tr}. \ep

\bcor \label{trans}
\begin{enumr}
\item $\ul \Con_n(\ul \Con_n(\gs))\vdash_\PA \ul \Con_n(\gs)$, for all $n<\gw$;
\item $\ul{\Con}_{\gw}(\ul{\Con}_\gw(\gs))\vdash_{\PA} \ul{\Con}_\gw(\gs)$.
\end{enumr}
\ecor

\bp\ Since the theories numerated by $\ul \Con_n(\gs)$ and $\ul
\Con_n(\ul \Con_n(\gs))$ are finite extensions of $\PA$, for a proof
of Statement (i) it is sufficient to show \beq \label{fin}
\PA+\Con_n(\ul \Con_n(\gs))\vdash \Con_n(\gs). \eeq Since
$\Con_n(\gs)$ is a $\Pi_{n+1}$-sentence, we can take in Lemma
\ref{reflc} $\phi=\Con_n(\gs)$ and $S=\PA+\phi$. This yields
statement \refeq{fin}.

For a proof of (ii), we show an informal version of this statement
by an argument formalizable in $\PA$. We must prove that, for each
$n<\gw$,
$$\PA+\Con_\gw(\ul \Con_\gw(\gs))\vdash \Con_n(\gs).$$ Using the monotonicity
and Statement (i) we reason as follows:
$$\ul{\Con}_{\gw}(\ul{\Con}_\gw(\gs))\vdash_{\PA} \ul{\Con}_n(\ul{\Con}_\gw(\gs)) \vdash_\PA
\ul{\Con}_n(\ul{\Con}_n(\gs))\vdash_{\PA} \ul{\Con}_n(\gs).$$ This
shows the claim. \ep

Corollary \ref{trans} shows the soundness of the third group of
rules of $\RCo$. As we mentioned above, the fourth group is actually
derivable from the first three and the monotonicity, so we can skip
it. We show the soundness of Axiom 5.

\bl If $n>m$ then $\PA\vdash \Con_n(\gs)\land \Con_m(\tau)\to
\Con_n(\gs\lor \ul \Con_m(\tau))$. \el

\bp\ We reason in $\PA$ as follows: If $\phi\in \Pi_{n+1}$ and
$\Box_{\gs\lor \ul \Con_m(\tau)}(\phi)$, then by the formalized
deduction theorem $\Box_\gs(\Con_m(\tau)\to \phi)$. Since $m<n$, the
formula $\Con_m(\tau)\to \phi$ belongs to $\Pi_{n+1}$. By
$\Con_n(\gs)$ we obtain $T_n(\gn{\Con_m(\tau)\to \phi})$ whence
$T_n(\gn{\Con_m(\tau)})\to T_n(\gn{\phi})$. Since $\Con_m(\tau)\in
\Pi_{n+1}$, from $\Con_m(\tau)$ we infer $T_n(\gn{\Con_m(\tau)})$.
Hence $T_n(\gn{\phi})$, as required. \ep

\bcor \label{j-ax} $\ul \Con_\gw(\gs)\lor \ul \Con_m(\tau)
\vdash_\PA \ul \Con_\gw(\gs \lor \ul \Con_m(\tau))$. \ecor

\bp\ Informally, we must prove, for each $n$, that
$$\PA+\Con_\gw(\gs)+ \Con_m(\tau)
\vdash  \Con_n(\gs \lor \ul \Con_m(\tau)).$$ We can assume $n>m$ and
then use the previous lemma. This argument is formalizable in $\PA$.
\ep

\bcor $\ul{\Con}_\gw(\gs)\vdash_\PA \gs$. \ecor

\bp\ We reason as follows:
\begin{eqnarray*}
\PA\vdash\ \Box_\gs(x) & \to &  \ex{n} (x\in \Pi_{n+1} \land \Box_\gs(x))\\
%& \to & \ex{n<x}\Box_{\Con_n(\gs)}(\Con_n(x)) \\
& \to & \ex{n} \Box_{\ul{\Con}_n(\gs)}(x) \\
& \to & \Box_{\ul{\Con}_\gw(\gs)}(x).\quad \boxtimes
\end{eqnarray*}
%\bigskip
This shows the soundness of the remaining Axiom 7 of $\RC$ and
completes the proof of Proposition \ref{sound}. \ep

\section{Kripke models for $\RCo$}

Kripke frames and models are understood in this paper in the usual
sense. A \emph{Kripke frame} $\cW$ for the language $\cL_S$ consists
of a non-empty set $W$ equipped with a family of binary relations
$(R_\ga)_{\ga\in S}$. A Kripke frame $\cW$ is called \emph{finite}
if so is $W$ and all but finitely many relations $R_\ga$ are empty.

A \emph{Kripke model} $\cW$ is a Kripke frame together with a
\emph{valuation} $v: W\times \Var \to \{0,1\}$ assigning a truth
value to each propositional variable at every node of $\cW$. As
usual, we write $\cW,x\fc A$ to denote that a formula $A$ is true at
a node $x$ of a model $\cW$. This relation is inductively defined as
follows:

\bi
\item $\cW,x\fc p \iff v(x,p)=1$, for each $p\in\Var$;
\item $\cW,x\fc \top$; \quad $\cW,x\fc A\land B \iff (\cW,x\fc
A \text{ and } \cW,x\fc B)$;
\item $\cW,x\fc \ga A \iff \ex{y}(xR_\ga y\text{ and }
\cW,y\fc A)$. \ei

We call a \emph{RJ$_S$-frame} a Kripke frame satisfying the
following conditions, for all $\ga,\gb\in S$ and all $x,y,z\in W$:
\bi
\item $xR_\ga yR_\gb z$ implies $x R_\gy z$, if $\gy=\min(\ga,\gb)$;
\hfill (polytransitivity)
\item $xR_\ga y$ and $xR_\gb z$ implies $y R_\gb z$, if $\ga>\gb$.
\hfill (condition J) \ei These conditions can be more succinctly
written as $R_\ga R_\gb\subseteq R_{\min(\ga,\gb)}$ and
$R_\ga^{-1}R_\gb\subseteq R_\gb$. An \emph{RC$_S$-frame} is an
RJ$_S$-frame that is \emph{monotone}, that is, $R_\ga\subseteq
R_\gb$, for each $\gb<\ga$. An \emph{RJ$_S$-model (RC$_S$-model)},
respectively, is a Kripke model based on an RJ$_S$-frame
(RC$_S$-frame).  We speak about RJ- and RC-frames and models
whenever $S=\gw+1$.

The persistence axiom $\gw A\vdash A$ does not correspond to a frame
condition.\footnote{Notice that the more familiar form of this axiom
is the principle $A\to \Box A$ which has no non-discrete frames.} We
call a Kripke model \emph{(downwards) persistent} if, for each
variable $p$,
$$\cW,x\fc p \text{ and } y R_\gw x \Imp \cW,y\fc p.$$
By a straightforward induction we obtain the following lemma. \bl
Let $\cW$ be a persistent Kripke model based on a polytransitive
frame. Then, for each positive formula $A$,
$$\cW,x\fc A \text{ and } y R_\gw x \Imp \cW,y\fc A.$$
\el

We say that a sequent $A\vdash B$ is \emph{true} in a Kripke model
$\cW$, if
$$\al{x\in \cW} (\cW,x\fc A\ \Imp\ \cW,x\fc B).$$
A logic $L$ is \emph{sound} for a class $\cC$ of Kripke models (of
the same signature), if every sequent $A\vdash B$ provable in $L$ is
true in any model from $\cC$. It is easy to see that our logics are
sound for their respective classes of models.

\bl \label{sound-Kr} \benr \item $\RJ_S$ is sound for the class of
all RJ$_S$-models;
\item $\RC_S$ is sound for the class of all RC$_S$-models;
\item $\RCo_S$ is sound for the class of all persistent RC$_S$-models.
 \eenr
\el A proof of this lemma is routine.

Notice that the frame conditions for the logics $\RJ_S$ and $\RC_S$
(that is, polytransitivity, condition J, and monotonicity) are
closure conditions. Therefore, for any Kripke frame
$\cW=(W,(R_\ga)_{\ga\in S})$ there is an RJ$_S$-frame (RC$_S$-frame)
$\ol\cW=(W,(\ol R_\ga)_{\ga\in S})$ such that \bi
\item $R_\ga\subseteq \ol R_\ga$, for each $\ga\in S$;
\item For any other RJ$_S$-frame (RC$_S$-frame) $(W,(R'_\ga)_{\ga\in S})$
with $R_\ga\subseteq R'_\ga$, for all $\ga\in S$, we have $\ol
R_\ga\subseteq R_\ga'$, for all $\ga\in S$. \ei The frame $\ol\cW$
is unique up to isomorphism. We call it the \emph{RJ$_S$-closure
(RC$_S$-closure)} of $\cW$.

\bex \label{ex} Consider a Kripke frame $\cW=(W,(R_\ga)_{\ga\leq
\gw})$ with $W=\{0,1,2\}$. Relation $R_\gw$ consists of two pairs $0
R_\gw 1$ and $0 R_\gw 2$, and the other relations are empty. Let
$\ol\cW$ be the RC-closure of $\cW$. It is easy to see that $\ol
R_\gw= R_\gw$, whereas, for each $n<\gw$, $\ol R_n\rst\{1,2\}$ is a
total relation, $0 \ol R_n 1$ and $0 \ol R_n 2$.

Further, we define $v(x,p)=0$ iff $x= 2$, and $v(x,q)=0$ iff $x=1$.
This makes $(\ol\cW,v)$ a downwards persistent RC-model falsifying
$\gw p\land \gw q \vdash \gw(p\land q)$ at $0$. By
Lemma~\ref{sound-Kr} we conclude $\gw p\land \gw q \nvdash_\RCo
\gw(p\land q)$. \eex

The completeness proofs in all these cases are also easy. As in
Dashkov~\cite{Das12}, we present an argument based on a (simplified)
version of filtrated canonical model.

Let $\Phi$ be a set of $\cL$-formulas. Denote
$\ell(\Phi):=\{\ga\leq\gw: \text{$\ga$ occurs in some
$A\in\Phi$}\}$. A set $\Phi$ is called \emph{adequate} if $\Phi$ is
closed under subformulas, $\top\in\Phi$ and \bi
\item If $\gb A\in\Phi$ and $\gb<\ga\in \ell(\Phi)$, then
$\ga A\in\Phi$;
\item For any variable $p$, if $p\in \Phi$ then $\gw p\in\Phi$.
\ei It is easy to see that any finite set of formulas can be
extended to a finite adequate set.

Let $\Gamma$ be a set of $\cL$-formulas and $L$ a logic. We shall
take for $L$ one of $\RJ$, $\RC$ or $\RCo$, or their restricted
versions in the language $\cL_S$ where $S\subseteq\gw+1$. We write
$\Gamma\vdash_L B$ if there are formulas $A_1,\dots, A_n\in \Gamma$
such that the sequent $A_1\land\cdots\land A_n\vdash B$ is provable
in $L$.

Fix an adequate set of formulas $\Phi$. An \emph{$L$-theory in
$\Phi$} is a set $x \subseteq \Phi$ such that $x \vdash_L A$ and
$A\in \Phi$ implies $A\in x $. Define a model $\cW_L/\Phi$ as
follows. The set of nodes $W_L/\Phi$ is the set of all $L$-theories
in $\Phi$.\footnote{In positive logic there is no harm in allowing
an `inconsistent' theory $x =\Phi$ as a node.} We stipulate that $x
R_\ga y$ iff $\ga\in\ell(\Phi)$ and the following conditions hold
for each formula $A$: \ben \renewcommand{\labelenumi}{R\theenumi.}
\item $A\in y$ and $\ga A\in \Phi$ implies $\ga A\in  x $;
\item $\gb A\in  y$ and $\ga A\in \Phi$
implies $\min(\ga,\gb)A\in x $;
\item $\gb<\ga$ and $\gb A\in  x $ implies $\gb A\in y$.
\een We also let $\cW_L/\Phi, x \fc p$ iff $p\in  x $, for any
$L$-theory $ x $.

\bl Suppose $L$ contains $\RJ_S$ with $S=\ell(\Phi)$. Then
$\cW_L/\Phi$ is an RJ$_S$-model. \el

\bp\ To check the polytransitivity assume $ x  R_\ga  y R_\gb z$ and
$\ga,\gb\in \ell(\Phi)$. We show $ x R_{\min(\ga,\gb)}  z$ by
checking R1--R3. If $A\in
 z$ and $\min(\ga,\gb)A\in\Phi$, then by the adequacy $\gb
A\in \Phi$ and hence $\gb A\in  y$. It follows that
$\min(\ga,\gb)A\in  x $.

For R2 notice that $\min(\gy,\ga,\gb)=\min(\gy,\min(\ga,\gb))$. If
$\gy A\in  z$ and $\min(\gy,\ga,\gb)A\in\Phi$ then by the adequacy
$\min(\gy,\gb) A\in \Phi$ and hence $\min(\gy,\gb) A\in  y$. This in
turn implies $\min(\gy,\gb,\ga) A\in  x $. Condition R3 is obviously
satisfied, as all three theories have the same formulas of the form
$\gb A$ for $\gb<\ga$.

Second, we check condition (J). Assume $ x R_\ga  y$ and $ x R_\gb
z$ with $\ga<\gb$. We show $ z R_\ga  y$. R1: If $A\in y$ and $\ga
A\in \Phi$ then $\ga A\in  x$. Since $\ga<\gb$ this implies $\ga
A\in z$. R2: If $\gy A\in  y$ and $\ga A\in \Phi$ then
$\min(\gy,\ga)A\in x$ whence $\min(\gy,\ga)A\in  z$ for the same
reason. R3 is, again, obvious. \ep

\bl \label{canon} For any $A\in\Phi$, $\cW_L/\Phi, x \fc A$ iff
$A\in x $. \el

\bp\ Induction on the build-up of $A$. If $A$ is a variable, $\top$
or has the form $B\land C$, the argument is obvious. Assume $A=\ga
B$.

If $ x \fc \ga B$ then, for some $ y$ such that $ x  R_\ga
 y$, we have $ y\fc B$. By IH it follows that $B\in
 y$ and hence $\ga B\in  x $.

Now assume $\ga B\in  x $. Let $\Delta:=\{\gb C : \gb C\in  x , \gb
<\ga\}$ and let $y$ be the deductive closure of $\Delta\cup \{B\}$
in $\Phi$. By the IH we have $ y\fc B$. We claim that $ x  R_\ga
 y$ which completes the argument.

Assume $D\in  y$, then $\Sigma, B\vdash_L D$ for some finite
$\Sigma\subseteq\Delta$. Then $\bigwedge \Sigma\land \ga B\vdash_L
\ga(\bigwedge\Sigma\land B)\vdash_L \ga D$. Hence, if $\ga D\in
\Phi$ then $\ga D\in  y$. Similarly, if $\gy D\in  y$ then $\Sigma,
B\vdash_L \gy D$. Then $\bigwedge \Sigma\land \ga B\vdash_L
\ga(\bigwedge\Sigma\land B)\vdash_L \ga\gy D\vdash_L
\min(\gy,\ga)D$. If $\ga D\in\Phi$ then $\min(\gy,\ga)D\in \Phi$
whence $\min(\gy,\ga)D\in x $. Finally, if $\gb<\ga$ and $\gb D\in x
$ then $\gb D\in\Delta$, hence $\gb D\in y$. \ep

\bl \benr
\item If $L$ contains the
monotonicity axiom and $S=\ell(\Phi)$, then $\cW_L/\Phi$ is an
RC$_S$-frame;
\item
If $L$ contains the persistence axiom, then $\cW_L/\Phi$ is
persistent. \eenr \el

\bp\ (i) Assume $ x  R_\ga  y$ and $\gb<\ga\in \ell(\Phi)$. We show
$x R_\gb y$ by checking the three conditions. If $A\in  y$ and $\gb
A\in \Phi$ then $\ga A\in \Phi$ by the adequacy of $\Phi$. Hence,
$\ga A\in  x $ and therefore by the monotonicity axioms $ x \vdash
\gb A$. Since $\gb A\in\Phi$ we obtain $\gb A\in x $. Similarly, if
$\gy A\in  y$ and $\gb A\in\Phi$ then $\ga A\in\Phi$ by the
adequacy. Therefore, $\min(\ga,\gy)A\in x $, whence by the
monotonicity axioms $ x \vdash\min(\gb,\gy)A$. Since both $\gy A$
and $\gb A$ are in $\Phi$, it follows that $\min(\gb,\gy)A\in x $
which proves the second condition. The third condition is obviously
satisfied.

(ii) Assume $ x  R_\gw  y$. If $ y\fc p$ then $p\in y$; by the
adequacy $\gw p\in\Phi$ and hence $\gw p\in  x $. It follows that $
x \vdash_L \gw p\vdash_L p$ and $ x \fc p$. \ep

Taking $\Phi=\cL$ and $\cW_L:= \cW_L/\Phi$ we obtain the
completeness of $\RJ$, $\RC$ and $\RCo$ w.r.t.\ their respective
classes of models.

\bt \label{comp} \benr
\item $A\vdash_\RJ B$ iff $A\vdash B$ is true in all RJ-models;
\item $A\vdash_\RC B$ iff $A\vdash B$ is true in all RC-models;
\item $A\vdash_\RCo B$ iff $A\vdash B$ is true in all persistent
RC-models. \eenr \et

\bp\ The three systems are sound by Lemma \ref{sound-Kr}. The
completeness is proved by observing that $\cW_L$, for each of the
three logics $L$, is a model of the corresponding type. Assume
$A\nvdash_L B$. Then letting $x$ denote the $L$-theory generated by
$A$ we have $B\notin x$, hence by Lemma \ref{canon} $\cW_L,x\nfc B$.
\ep

Next we discuss the finite model property of the three logics. For
$\RJ$ the answer is obvious, but for $\RC$ and $\RC\gw$ we have a
small complication due to the fact that modality $\gw$ is present in
the language.

\bcor $A\vdash_\RJ B$ iff $A\vdash B$ is true in all finite
RJ-models. \ecor

\bp\ Assume $A\nvdash_\RJ B$, let $\Phi$ be a finite adequate set of
formulas containing both $A$ and $B$. We have $\cW_\RJ/\Phi,x\fc A$
and $\cW_\RJ/\Phi,x\nfc B$. Moreover, $\cW_\RJ/\Phi$ is a finite
RJ$_S$-model, where $S=\ell(\Phi)$. By putting $R_\ga:=\emptyset$,
for any $\ga\notin\ell(\Phi)$, we expand $\cW_\RJ/\Phi$ to an
RJ-model in $\cL$ falsifying $A\vdash B$. \ep

A similar argument does not quite work for $\RC$, as the expansion
by empty relations leads, in general, outside the class of
RC-models. However, for a finite signature $S$ we do have an analog
of Theorem \ref{comp}.

\bcor \label{fincomp} Suppose $S\subseteq\gw+1$ is finite. \benr
\item $A\vdash_{\RC_S} B$ iff $A\vdash B$ is true in all finite
RC$_S$-models; \item $A\vdash_{\RC\gw_S} B$ iff $A\vdash B$ is true
in all finite persistent RC$_S$-models.\eenr \ecor

\bl \label{expand} Let $S\subseteq\gw+1$ and $\ga\notin S$. Any
RC$_S$-model can be expanded to an RC$_{S\cup \{\ga\}}$-model. \el

\bp\ Let $\cW=(W,(R_\gb)_{\gb\in S})$ be a given RC$_S$-model.
Denote: $S^+:=\{\gb\in S:\ga<\gb\}$ and $S^-:=\{\gb\in S:\gb<\ga\}$.
If $S^+=\emptyset$ we can put $R_\ga:=\emptyset$. Otherwise, for any
relation $R$ on $W$, denote $$R':= R\cup RR\cup \bigcup_{\gb\in
S^+}R_\gb^{-1}R.$$ Further, define $R_\ga:=\bigcup_{n\in \gw}
R_\ga^n$, where $R_\ga^0:= \bigcup_{\gb\in S^+} R_\gb$;
$R_\ga^{n+1}:=(R_\ga^n)'.$

Notice that $R\subseteq R'$, for any $R$. It follows that
$R_\gb\subseteq R_\ga^0\subseteq R_\ga$, for each $\gb\in S^+$. By
the construction, $R_\ga$ is transitive and
$R_\gb^{-1}R_\ga\subseteq R_\ga$, hence condition (J) is satisfied
for all $\ga<\gb\in S^+$. Moreover, the polytransitivity follows
from the transitivity and the monotonicity properties. Therefore,
$(W,(R_\gb)_{\gb\in S^+\cup\{\ga\}})$ is an
RC$_{S^+\cup\{\ga\}}$-model.

To complete the argument we have to show that $(W,(R_\gy)_{\gy\in
S^-\cup\{\ga\}})$ is an RC$_{S^-\cup\{\ga\}}$-model. To this end we
prove that, for each $n$ and $\gy\in S^-$, \ben
\item $R_\ga^n\subseteq R_\gy$;
\item $(R_\ga^n)^{-1} R_\gy\subseteq R_\gy$.
\een Both statements are verified by induction on $n$. The basis of
induction holds, since the original model was an RC$_S$-model.
Assume the statements hold for $R=R_\ga^n$ and consider
$R'=R_\ga^{n+1}$.

1. We have $R\subseteq R_\gy$ by the IH. Further, $RR\subseteq R_\gy
R_\gy\subseteq R_\gy$, since $R_\gy$ is transitive. For any $\gb\in
S^+$, $R_\gb^{-1} R\subseteq R_\gb^{-1}R_\gy\subseteq R_\gy$, since
condition (J) holds in $\cW$. Hence, $R'=R\cup RR\cup
\bigcup_{\gb\in S^+}R_\gb^{-1}R \subseteq R_\gy$, as required.

2. We have $R^{-1}R_\gy\subseteq R_\gy$ by the IH. Further,
$(RR)^{-1}R_\gy=R^{-1}(R^{-1}R_\gy)\subseteq R^{-1}R_\gy\subseteq
R_\gy$. Finally, for any $\gb\in S^+$, $(R_\gb^{-1}R)^{-1}
R_\gy=R^{-1}R_\gb R_\gy\subseteq R^{-1}R_\gy\subseteq R_\gy$.
Therefore, $(R')^{-1} R_\gy\subseteq R_\gy$, as required. \ep

\brem The given proof also works for the more general analogs of
$\RC$, e.g., for logics with linearly ordered sets of modalities
(see \cite{BFJ12}). \erem

Taking into account that expansions of persistent models are
persistent, we obtain the following theorem for both $\RC$ and
$\RC\gw$.

\bt \label{conserv} Let $L$ be either $\RC$ or $\RC\gw$. The
following statements are equivalent: \benr
\item $A\vdash_L B$;
\item $A\vdash_{L_U} B$, for some finite $U\subseteq \gw+1$;
%\item $A\vdash_{\RC_M} B$ where $M:= \max(S)$;
\item $A\vdash_{L_S} B$ where $S=
\ell(\{A,B\})$. \eenr \et

\bp\ Clearly, (iii) implies (i), and (i) implies (ii) since a finite
derivation may only contain finitely many different modalities. We
prove that (ii) implies (iii). Assume $A\nvdash_{L_S} B$. By
Corollary \ref{fincomp} there is a finite RC$_S$-model $\cW$
falsifying $A\vdash B$ (which is persistent if $L=\RC\gw$). Assume
any finite $U$ be given. We may assume $S\subseteq U$ (otherwise
clearly $A\nvdash_{L_U} B$). By Lemma \ref{expand}, $\cW$ can be
expanded to an RC$_U$-model falsifying the same sequent. Hence,
$A\nvdash_{L_U} B$. \ep

%\bcor $\RC$ is conservative over $\RC_S$, for any $S\subseteq\gw+1$.
%\ecor

Thus, even though we do not have the finite model property for $\RC$
and $\RC\gw$ in the full language, these logics are conservatively
approximated by their fragments with this property. Together with
Corollary \ref{fincomp} this yields

\bcor The systems $\RC$ and $\RCo$ are decidable. \ecor

For the logics $\RJ$ and $\RC$ a sharper result can be stated. As we
have seen, the question whether a sequent $A\vdash B$ is provable in
such a logic $L$ is equivalent to the same question for the logic
$L_S$ with $S=\ell(\{A,B\})$. However, for any finite $S$, the logic
$L_S$ is modulo renaming of modalities the same logic as $L_n$ for
$n=|S|$ (we identify $n$ with the set $\{0,\dots,n-1\}$). The
systems $L_n$ are shown to be polytime decidable \cite{Das12}.
Therefore, we obtain

\bcor The systems $\RJ$ and $\RC$ are polytime decidable. \ecor

The same result holds for $\RCo$, however we cannot directly refer
to Dashkov's theorem. This question is considered in the next
section, where we also obtain somewhat sharper complexity estimates
for the cases $\RJ$ and $\RC$. The material of that section, up to
Theorem \ref{cantree}, is due to Dashkov \cite{Das12}.

\section{Polytime decidability of $\RCo$}
We have to develop some combinatorial techniques to deal with
positive logics. It allows one to state the Kripke completeness
results in a sharper form, from which the complexity bounds are
easily read off.

Let $\cW$ be a Kripke model and $a\in W$. The submodel $\cW_a$ of
$\cW$ generated by $a$ is obtained by restricting all the relations
and the valuation of $\cW$ to the set of all nodes $x\in W$ such
that there is a path $a=x_0Rx_1R\dots Rx_n=x$ where
$R=\bigcup_{\ga\in S} R_\ga$. A model $\cW$ is called \emph{rooted}
if it has a distinguished element $a$ (called the \emph{root}) such
that $\cW_a=\cW$. We notice that in polytransitive rooted frames
every node is reachable from the root in one step.

\bd We can associate with each positive formula $A$ a rooted
treelike Kripke model $T[A]$ in the signature $\ell(A)$ called its
\emph{canonical tree}. It is essentially the parse tree of $A$
viewed as a Kripke model.

If $A$ is a variable or $\top$, then $T[A]$ is a one-point model
$\{a\}$ with the empty relations, and the only variable true at $a$
is $A$.

If $ A=B\land  C$ then $T[A]$ is obtained from the disjoint union of
the models $T[B]$ and $T[C]$ by identifying the roots. We declare
any variable $p$ true at the root of $T[A]$ iff it is true at the
root of either $T[B]$ or $T[C]$.

If $A=\ga B$ then $T[A]$ is obtained from $T[B]$ by adding a new
root $r$ (where all variables are false), from which the root of of
$T[B]$ is $R_\ga$-accessible. \ed

We write $T[A]\fc \phi$ if $\phi$ is true at the root of $T[A]$.
Then, one can easily verify the following properties: \bi
\item Each $R_\ga$ on $T[A]$ is an irreflexive forest-like binary relation;
\item $T[A]\fc  A$.
\ei

\bd A \emph{homomorphism} of a Kripke model $\cV$ into a Kripke
model $\cW$ (of the same signature $S$) is a function $f:V\to W$
such that \bi
\item $\al{x,y\in V} (xR_\ga y \Imp f(x)R_\ga f(y))$, for each
$\ga\in S$;
\item If $\cV,x\fc p$ then $\cW,f(x)\fc p$, for each variable
$p$. \ei

Let $\cV$ and $\cW$ be rooted Kripke models. A \emph{simulation} of
$\cV$ by $\cW$ is a homomorphism $f:\cV\to \cW$ mapping the root of
$\cV$ to the root of $\cW$. \ed

\bl \label{mon} If $A$ is strictly positive and $f$ is a
homomorphism of $\cV$ into $\cW$, then $$\cV,x\fc A \Imp \cW,f(x)\fc
A.$$ \el

\bl \label{sim} $\cW,x\fc B$ iff there is a homomorphism $f:T[B]\to
\cW$ mapping the root of $T[B]$ to $x$. \el

\bp\ Suppose $f:T[B]\to \cW$ is such a homomorphism. We have $T[B],
r\fc B$ where $r$ is the root of $T[B]$. Since $B$ is strictly
positive and $f(r)=x$, by Lemma \ref{mon}, $\cW,x\fc B$.

Suppose $\cW,x\fc B$. We construct a homomorphism $f:T[B]\to \cW$ by
induction on the complexity of $B$. If $B$ is a variable or $\top$,
the claim is obvious.

If $B= C\land  D$ then $\cW,x\fc C, D$. By the IH, there are
homomorphisms $f, g$ of the models $T[C]$ and $T[D]$ into $\cW$
mapping their respective roots to $x$. The homomorphism of $T[B]$
maps its root to $x$ and is defined as the union of $f$ and $g$
everywhere else on $T[B]$. We note that if $T[B]\fc p$ then either
$T[C]\fc p$ or $T[D]\fc p$, by the definition of $T[B]$. In either
case we have $\cW,x\fc p$, therefore the variable condition at the
root is met and we have a homomorphism of $T[B]$ into $\cW$.

If $ B=\ga C$ and $\cW,x\fc  B$, then there is a node $y\in \cW$
such that $x R_\ga y$ and $\cW,y\fc C$. By the IH, there is a
homomorphism of $T[C]$ into $\cW$ mapping its root to $y$. We extend
it by mapping the root of $T[B]$ to $x$. All the variables are false
at the root of $T[B]$, so the variable condition is met. \ep

Let $\Rc_S[A]$ ($\rj_S[A]$) denote the RC$_S$-closure (respectively,
RJ$_S$-closure) of $T[A]$, where $S\supseteq \ell(A)$.

\bt \label{cantree} \benr
\item $A\vdash_\RJ B$ iff $\rj_S[A]\fc B$, where $S=\ell(A)$;
\item $A\vdash_\RC B$ iff $\Rc_S[A]\fc B$, where $S=\ell(\{A,B\})$.
\eenr \et

\bp\ We prove Statement (ii). The case of $\RJ$ is similar but
simpler.

(only if) Since $T[A]\fc A$ and the relations of $\Rc_S[A]$ extend
those of $T[A]$, we have $\Rc_S[A]\fc A$. By Theorem \ref{conserv},
$A\vdash_\RC B$ implies $A\vdash_{\RC_S}B$. Hence, by Corollary
\ref{fincomp}, $\Rc_S[A]\fc B$.

(if) Assume $A\nvdash_\RC B$. There is a rooted RC$_S$-model $\cW$
such that $\cW\fc A$ and $\cW\nfc B$. By Lemma \ref{sim}, there is a
simulation $f: T[A]\to \cW$. Since $\cW$ is an RC$_S$-model, $f$
lifts to a simulation of $\Rc_S[A]$ by $\cW$. In fact, we can define
on $T[A]$ new relations $R_\ga'$ by letting $x R_\ga' y$ iff $f(x)
R_\ga f(y)$ in $\cW$. Then $(T[A],(R_\ga')_{\ga\in S})$ is an
RC$_S$-model with $R_\ga\subseteq R_\ga'$, for all $\ga\in S$.
Hence, denoting by $R_\ga''$ the relations of $\Rc_S[A]$, we obtain
$R''_\ga\subseteq R_\ga'$, for each $\ga\in S$. It follows that
$\cW$ simulates $\Rc_S[A]$ by $f$. Then, since $\cW\nfc B$, we
conclude that $\Rc_S[A]\nfc B$, by Lemma \ref{mon}. \ep

\brem The proof of Theorem \ref{cantree} provides an alternative way
of showing the finite model property for the logics $\RC_S$ and
$\RJ_S$. \erem

Theorem \ref{cantree} yields an efficient decision procedure for the
logics $\RC$ and $\RJ$. Firstly, given a positive formula $A$ we let
$S=\ell(\{A,B\})$ and build the model $\Rc_S[A]$. Secondly, we check
if $B$ is satisfied at the root of this model. To estimate the
complexity of this procedure we need to be more specific about the
chosen computation model.

We consider random access machines (see \cite{CoRe77}) and assume
that any register can hold (the code of) any symbol including the
variables and the modalities. To simplify the estimates we count the
size of any symbol as one, and we assume that the elementary
operations such as reading and writing a symbol, as well as the
comparison of symbols, cost a constant amount of time. We are going
to estimate the number of elementary steps needed to decide whether
$A\vdash_{\RC} B$. (Representing the variables and the modalities
more faithfully would introduce a logarithmic factor into our
estimates.) First, we estimate the time needed to build the model
$\Rc_S[A]$ given $A$.

We support a data structure for a positive formula $A$ (and for the
corresponding Kripke model $T[A]$) with the arrows represented by
pointers. The arrows are labeled by the elements of $S$, the nodes
are labeled by the variables of $A$. We can also realize these
labels as pointers to some extra nodes representing the variables
and the modalities, respectively. We assume that there is a fixed
ordering of arrows outgoing from any given node of the tree $T[A]$
(which respects the left-to-right ordering of the corresponding
subformulas of $A$). It is well-known that we can very efficiently
(in a linear number of steps) parse the formula $A$ to build such a
tree.

Next we bring $T[A]$ to a special \emph{ordered} form. Let
$\cL_{\geq m}$ denote the language $\cL_U$ with $U=[m,\gw]$. A
formula will be called a \emph{fact} if it is either $\top$ or a
conjunction of variables. Ordered formulas are defined inductively.

\bd A formula $A$ is \emph{ordered} if it has the form $A=F\land
\bigwedge_{i< k} m_i A_i$ for some $k$ (assuming $A=F$ if $k=0$),
where

\benr
\item $F$ is a fact;
\item For each $i$, $A_i\in \cL_{\geq m_i}$ and $A_i$ is ordered;
\item $m_0 \geq m_1\geq  \dots \geq  m_{k-1}$.
\eenr \ed

\bl \label{or-eq} Every positive formula $A$ is $\RJ$-equivalent to
an ordered one.\el

\bp\ Induction on the build-up of $A$. The basis of induction and
the case of conjunction are easy. Suppose $A=m B$. By the induction
hypothesis we may assume $B$ ordered, that is, $B=F\land
\bigwedge_{i< k} m_i B_i$. If $mB$ is not ordered, there is an $i<k$
such that $m_i<m$. Let $s$ be the minimal such $i$. Then $mB$ is
equivalent to $\top \land m(F\land \bigwedge_{i< s} m_i B_i)\land
\bigwedge_{i=s}^{k-1}  m_i B_i,$ which is ordered. \ep

We notice that if an ordered formula $B$ is obtained from $A$ by the
recursive procedure described in Lemma \ref{or-eq}, then the number
of nodes in $T[B]$ is the same as in $T[A]$. One can also easily
prove that $\Rc_S[B]$ will, in fact, be isomorphic to $\Rc_S[A]$.

The algorithm of ordering a formula is similar to that of sorting a
string, and it is easy to obtain a rough quadratic upper bound, a
detailed proof of which we omit.

\bl \label{order} Any formula $A$ can be ordered in $O(|A|^2)$
steps. \el

An ordered formula $A$ can be written in the following form: \beq
A=F\land \bigwedge_{i<k}\bigwedge_{j<n_i} m_i A_{ij}, \label{fm}\eeq
with $m_0>m_1>\dots >m_{k-1}$, $A_{ij}\in\cL_{\geq m_i}$ ordered and
$F$ a fact. Then $\Rc_S[A]$ can be characterized as follows.

\bl  \label{closure} If $A$ of the form \refeq{fm} is ordered then
$\Rc_S[A]$ consists of the disjoint union of the models
$\Rc_{S_i}[A_{ij}]$, for all $i< k$ and $j<n_i$, augmented by a new
root $a$, where $S_i:=S\cap [m_i,\gw]$ and $a\fc F$. In addition to
all the relations inherited from the models $\Rc_{S_i}[A_{ij}]$,
only the following relations hold in $\Rc_S[A]$:

\ben
\item $a R_{n} x$, for each $i< k$, $n\leq m_i$, $j<n_i$ and $x\in \Rc_{S_i}[A_{ij}]$;
\item $x R_n y$, for each $i< k$, $m_{i+1}\leq n<m_i$, and $x,y\in
\bigcup_{p\leq i}\bigcup_{j< n_p} \Rc_{S_p}[A_{pj}]$ (where we
formally let $m_k=0$);
\item $x R_{m_i} y$, for each $i< k$, $y\in \bigcup_{j<n_i}\Rc_{S_i}[A_{ij}]$
and $x\in \bigcup_{p<i}\bigcup_{j<n_p} \Rc_{S_p}[A_{pj}]$. \een

\el

\bp\ It is easy to see that all the relations mentioned in items
1--3 must hold in $\Rc_S[A]$.

1. By the polytransitivity we have $aR_{m_i} x$, for each
$x\in\Rc_{S_i}[A_{ij}]$. Then, by the monotonicity, $aR_n x$, for
all $n\leq m_i$.

2. If $x,y\in\bigcup_{p\leq i}\bigcup_{j< n_p} \Rc_{S_p}[A_{pj}]$
then $aR_{m_i} x,y$ by Item 1, since $m_p\geq m_i$, for each $p\leq
i$. In particular, for each $n<m_i$, there holds $aR_n y$. Then by
property (J) we obtain $x R_n y$.

3. For any $x,y$ as specified we have $aR_{m_i} y$ and $a
R_{m_{i-1}} x$ by Item 1. Since $m_{i-1}>m_i$, by (J) we obtain $x
R_{m_i} y$.

It is also a routine but somewhat lengthy check that the model
described in Lemma \ref{closure} is, indeed, an RC$_S$-model. Hence,
it must coincide with $\Rc_S[A]$. \ep

A similar but much simpler characterization holds for $\rj_S[A]$. In
this case, we do not need to assume that $A$ is ordered. If $x,y\in
T[A]$ let $x\sqcap y$ denote the greatest lower bound of $x,y$, that
is, the unique node $z$ such that there are oriented paths from $z$
to $x$ and to $y$ without any shared edges. Given a nonempty path
$P$ let $m(P)$ denote the minimal modality label occurring on $P$.

\bl \label{rj-cl} Let $x,y\in T[A]$, $S=\ell(A)$, and let $X$ and
$Y$ be the uniquely defined paths from $x\sqcap y$ to the nodes $x$
and $y$, respectively. Then $x R_n y$ holds in $\rj_S[A]$ iff either
$x\sqcap y = x$, $x\neq y$ and $m(Y)=n$, or both $X$ and $Y$ are
nonempty and $m(X)>m(Y)=n$. \el

We omit a routine proof. Lemmas \ref{closure} and \ref{rj-cl} yield
efficient algorithms to build the models $\Rc_S[A]$ and $\rj_S[A]$.

\bl \label{clos} For any ordered formula $A$, the model $\Rc_S[A]$
can be constructed in $O(|A|^2\cdot |S|)$ many steps. \el

\bp\ In the course of constructing $\Rc_S[A]$ we add arrows to the
initial model $T[A]$ in a systematic way. Assume
$A=F\land\bigwedge_{i<k} \bigwedge_{j<n_i} m_iA_{ij}.$ Since
$A_{ij}\in\cL_{S_i}$, we can apply the procedure recursively to
build $\Rc_{S_i}[A_{ij}]$, for each $i< k$ and $j<n_i$. Then we join
these models by a common root and add the arrows according to
clauses 1--3 of Lemma~\ref{closure}. Apart from the computation time
needed to build the models $\Rc_{S_i}[A_{ij}]$, this requires only a
linear number of steps in the number of added arrows. Further,
notice that we never add an arrow twice (clauses 1--3 enumerate
distinct fresh arrows because of the choice of the sets $S_i$).
Thus, the total number of steps in the whole computation is linear
in the total number of arrows in the model $\Rc_S[A]$ which can be
roughly estimated by $O(|A|^2\cdot |S|)$. \ep

\bt \label{bound} The logics $\RJ$ and $\RC$ are decidable in time
bounded by a polynomial (of degree three and four, respectively) in
the length of the sequent $A\vdash B$. \et

\bp\ It is well-known that the problem of checking whether a modal
formula $\phi$ is true at the root of a finite Kripke model $\cW$ in
a finite signature $S$ is solvable in time $O(\|W\|\cdot |\phi|)$,
where $\|W\|$ denotes the sum of $|W|$ and the number of pairs
$(x,y)$ such that $xR_\ga y$, for some $\ga\in S$ (see
\cite[Proposition 3.1]{HM92}).\footnote{$\|W\|$ only measures the
complexity of the \emph{frame}, while the number of variables is
accounted for in $|\phi|$.}

Letting $S=\ell(\{A,B\})$ and $n=|S|$ we can estimate $|\Rc_S[A]|$
by $|A|$ and $\|\Rc_S[A]\|$ by $O(|A|^2\cdot n)$. This yields a
bound of the form $O(|A|^2\cdot n\cdot |B|)$ on the complexity of
checking whether $\Rc_S[A]\fc B$. By Lemmas \ref{order} and
\ref{clos} we have the same bound on the complexity of the original
problem $A\vdash_{\RC} B$. Noting that $n\leq |A|+|B|$ yields a
fourth degree polynomial bound in the length $|A|+|B|$ of the input.

For the logic $\RJ$ this can be slightly improved. By Lemma
\ref{rj-cl}, we can observe that in the graph $\rj_S[A]$ (where
$S=\ell(A)$) there is no more than one arrow between any pair of
points. This yields a bound $O(|A|^2)$ on $\|\rj_S[A]\|$ and on the
complexity of constructing this model. Consequently, the
derivability problem can be solved in $O(|A|^2\cdot |B|)$ many
steps. \ep

\brem Since the input of the problem is naturally divided into two
parts $A$ and $B$, measuring the complexity in terms of two
parameters $|A|$ and $|B|$ appears to be more meaningful than
expressing it in terms of the total length of the input. Thus, the
more informative  bounds are $O(|A|^2\cdot n\cdot |B|)$ for $\RC$
and $O(|A|^2\cdot |B|)$ for $\RJ$. \erem

Next we turn to the logic $\RCo$. We define $\Rc\gw_S[A]$ as the
model whose frame coincides with that of $\Rc_S[A]$ and whose
valuation function $v'$ satisfies:
$$v'(x,p)=1 \iff (v(x,p)=1 \text{ or }
\ex{y}(x R_\gw y \text{ and } v(y,p)=1)),$$ where $v$ is the
valuation of $\Rc_S[A]$. Clearly, $\Rc\gw_S[A]$ is persistent.

\bt $A\vdash_{\RCo} B$ iff $\Rc\gw_S[A]\fc B$, where
$S=\ell(\{A,B\})$. \et

\bp\ The proof is similar to that of Theorem \ref{cantree}, we
elaborate the (if) part.

If $A\nvdash_{\RCo} B$, there is a rooted persistent RC$_S$-model
$\cW$ such that $\cW\fc A$ and $\cW\nfc B$. By Lemma \ref{sim},
there is a simulation $f: T[A]\to \cW$. As before, since $\cW$ is an
RC$_S$-model, $f$ lifts to a simulation of $\Rc_S[A]$ by $\cW$. We
claim that $f$ also lifts to a simulation of $\Rc\gw_S[A]$ by $\cW$.
Assume $\Rc\gw_S[A]\fc p$. If $\Rc_S[A]\fc p$ then $\cW,f(x)\fc p$
and there is nothing to prove. If $x R_\gw y$ in $\Rc_S[A]$ and
$\Rc_S[A],y\fc p$, then $\cW,f(y)\fc p$. Hence, by the persistence
of $\cW$, we obtain $\cW,f(x)\fc p$ and the claim is proved.
Therefore, since $\cW\nfc B$, we have $\Rc\gw_S[A]\nfc B$, as
required. \ep

Now we notice that $\|\Rc\gw_S[A]\|$ has the same bound as
$\|\Rc_S[A]\|$. Hence, we obtain

\bcor The logic $\RC\gw$ is decidable in time bounded by a
polynomial (of degree four) in the length of the sequent $A\vdash
B$. \ecor

\section{Irreflexive models}
The Solovay construction works with the irreflexive models.
Therefore we would like to have a characterization of $\RC$ and
$\RC\gw$ in terms of suitable irreflexive models. We modify the
construction of the canonical model from the previous section. This
modification is similar to the one given by Dashkov \cite{Das12}
which in turn derives from the work of Japaridze \cite{Dzh86} and
Ignatiev \cite{Ign93}.

Let $L$ be a logic containing $\RC$ and let $\Phi$ be an adequate
set of formulas. We work in the setup of the previous section. We
define a Kripke model $\cW'_L/\Phi$ which coincides with
$\cW_L/\Phi$ but for the definition of the relations. We stipulate
that $x R'_\ga y$ in $\cW'_L/\Phi$ iff $x R_\ga y$ and following
condition holds: \ben
\renewcommand{\labelenumi}{R\theenumi.} \setcounter{enumi}{3}
\item There is a formula $\ga A\in x$ such that $\ga A\notin y$.
\een

The model $\cW_L'/\Phi$ has the following properties.

\bl \benr
\item $\cW_L'/\Phi$ is an RJ-model;
\item All $R'_\ga$ are irreflexive;
\item $R'_\ga=\emptyset$, for $\ga\notin \ell(\Phi)$.
%\item $\cW_L'/\Phi$ is persistent for $L=\RCo$.
\eenr \el

In order to prove the canonical model lemma we need an additional
fact.

\bl \label{irref} For any $A$ and $\ga$, $A\nvdash_{\RCo} \ga A$.
\el

\bp\ The simplest proof of this fact involves arithmetical
interpretations. Assume $A\vdash_{\RCo} \ga A$. Fix an arithmetical
interpretation $*$ mapping all variables to the standard numeration
of $\PA$. By Proposition \ref{sound} we obtain $A^*\vdash_\PA
\Con_\ga(A^*)$. Since $A$ is a positive formula, $A^*$ is (a
numeration of) a sound arithmetical theory $T$. However, this
contradicts G\"odel's second incompleteness theorem for $T$. \ep

\bl \label{canon1} For any $A\in\Phi$, $\cW'_L/\Phi, x \fc A$ iff
$A\in x $. \el

\bp\ The proof is similar to that of Lemma \ref{canon}. We only
treat somewhat differently the case $A=\ga B$, the `if' part.

Assume $\ga B\in  x $. As before let $\Delta:=\{\gb C : \gb C\in  x
, \gb <\ga\}$ and let $ y$ be the deductive closure of $\Delta\cup
\{B\}$ in $\Phi$. By the IH we have $ y\fc B$. We claim that $ x
R'_\ga y$ which completes the argument.

We already know from Lemma \ref{canon} that $x R_\ga y$. To check R4
it is sufficient to observe that $\ga B\in x$ but $\ga B\notin y$.
If $\ga B\in y$ we would obtain $\Sigma, B\vdash_L \ga B$, for some
finite $\Sigma\subseteq\Delta$, but then $\bigwedge\gS\land
B\vdash_L \bigwedge\gS\land \ga B\vdash_L \ga(\bigwedge \Sigma \land
B)$, contradicting Lemma \ref{irref}. \ep

A model $\cW$ is called \emph{$\Phi$-monotone}, if for any $\ga
A\in\Phi$ and $\gb\in\ell(\Phi)$ such that $\ga<\gb$, $\cW'_L/\Phi,x
\fc \gb A$ implies $\cW'_L/\Phi,x \fc \ga A$.

\bl \benr \item $\cW'_L/\Phi$ is $\Phi$-monotone;
\item $\cW_L'/\Phi$ is persistent if $L$ contains $\RCo$.
\eenr \el

\bp\ (i) Assume $\cW'_L/\Phi,x \fc \gb A$, $\ga A\in \Phi$ and
$\ga<\gb\in\ell(\Phi)$. By the adequacy of $\Phi$ we have $\gb
A\in\Phi$. Then by Lemma \ref{canon1} we obtain $\gb A\in x$. Hence,
$x\vdash_L \ga A$ and since $\ga A\in\Phi$ also $\ga A\in x$. This
yields $\cW'_L/\Phi,x \fc \ga A$ by Lemma \ref{canon1}.

Statement (ii) is obvious by Lemma \ref{canon1}, since $\cW_L/\Phi$
is persistent. \ep

We summarize the information obtained so far for $L=\RCo$.

\bpr \label{irmodel} Let $\Phi$ be a finite adequate set. Then there
is a finite model $\cW$ such that \benr
\item $\cW$ is an irreflexive  RJ-model; \item $R_\ga=\emptyset$, for all
$\ga\notin \ell(\Phi)$;
\item $\cW$ is \/$\Phi$-monotone and persistent;
\item For any $\RCo$-theory $\Gamma$ in $\Phi$ there is a node
$x\in \cW$ such that, for any formula $A$, $A\in\Gamma$ iff
$\cW,x\fc A.$
%\item $\cW,x\fc A$ and $\cW,x\nfc B$, for some $x\in \cW$;
\eenr \epr

\section{Arithmetical completeness}

\bt \label{arcomp} For any sequent $A\vdash B$ the following
statements are equivalent: \benr \item $A\vdash B$ is provable in
$\RCo$; \item $A^*\vdash_{\PA} B^*$, for all arithmetical
interpretations $*$;
\item $A^*\vdash B^*$, for all arithmetical interpretations $*$.
\eenr \et

\bp\ The implication from (i) to (ii) is Proposition \ref{sound}.
Statement (ii) trivially implies (iii). To infer (i) from (iii) we
argue by contraposition and assume $A\nvdash_\RCo B$. Consider a
finite adequate set $\Phi$ containing $A,B$, and let $\cW$ be a
finite Kripke model satisfying the conditions of Proposition
\ref{irmodel}. It falsifies $A\vdash B$ at some node $x$. We can
restrict $\cW$ to the submodel generated by $x$, so that $\cW$ is
rooted and falsifies $A\vdash B$ at the root.

\bigskip
Now we proceed to a Solovay-type construction. As usual, we identify
the nodes of $\cW$ with a finite set of natural numbers $\{1,\dots,
N\}$ so that $1$ is the root. We then attach a new root $0$ to $\cW$
by  stipulating that $0 R_0 x$, for each $x\in W$. The valuation of
variables at $0$ will be the same as in $1$, this ensures that the
new model is persistent. Abusing notation we denote this model by
the same letter $\cW$. We also assume that the $R_\ga$ relations and
the forcing relation $x\fc C$ on $\cW$ are arithmetized in a natural
way by bounded (even open) arithmetical formulas.

We fix an arithmetical formula $\Prf_n(x,y)$ naturally expressing
that \emph{$y$ is a proof of a formula $x$ from the axioms of
\,$\PA$ and true \/$\Pi_n$-sentences}. The formula $\Prf_n(x,y)$ has
logical complexity $\Delta_{n+1}$ in $\PA$. Without loss of
generality we may also assume that $\Prf_n$ is chosen in such a way
that each number $y$ is a proof of at most one formula, and that any
provable formula has arbitrarily long proofs. These properties are
also assumed to hold provably in $\PA$.

The formula $\Box_n(x):=\ex{ y}\Prf_n(x,y)$ expresses that $x$ is
provable in $\PA$ from the set of all true $\Pi_{n}$-sentences. We
usually write $\Box_n\phi$ for $\Box_n(\gn{\phi})$. It is easy to
see that $\Con_n(\gs_\PA)$ is equivalent to $\neg\Box_n\bot$.

\bd\label{Sol} Let $M$ denote the maximal modality $m<\gw$ occurring
in $\Phi$, if there is such an $m$, and $0$ otherwise. We define a
family of Solovay-style functions $h_n:\gw\to \cW$, for all $n<\gw$,
as follows: $h_n(0)=0$ and
$$
h_n(x+1)= \begin{cases}
y, & \text{if $h_i(x)\neq h_i(x+1)=y$, for some $i<n$; otherwise} \\
z, & \text{if $\ex{k\geq \max(M,n)}\Prf_n(\gn{\ell_k\neq
\bar{z}},x)$ and}
\\  & \text{either $h_n(x)R_n z$ or $h_n(x)R_\gw z$;} \\
h_n(x), & \text{otherwise}.
\end{cases}
$$
\ed Here $\ell_k$ denotes the limit of the function $h_k$. The
functions $h_n$ can be defined in such a way as to satisfy the
following conditions: \bi
\item The graph of each $h_n$ is definable by a
formula $H_n$ which is $\Delta_{n+1}$ in $\PA$;
\item The function $\phi: n\mapsto
\gn{H_n}$ is primitive recursive;
\item Each $h_n$ satisfies the clauses of Definition
\ref{Sol} provably in $\PA$. \ei

The definition of the functions $h_n$ can be arranged as a solution
of a fixed point equation in $\PA$ using the standard methods. The
details are given in the Appendix.

Informally, the behavior of the functions $h_n$ can be described as
follows. The functions with lower index have higher priority,
therefore whenever $h_m$ makes a move to $y$, all functions $h_n$
with $n>m$ do the same. Otherwise, $h_n$ moves like the usual
Solovay function, but for the following peculiarities: \bi
\item $h_n$ also reacts to proofs of the limit statements
for all functions of lower priority (not only to those of itself);
\item $h_n$ is not only allowed to move along the $R_n$ relation
but also along $R_\gw$. \ei

\bl \label{incr} For each $n,m$, provably in $\PA$,
\begin{enumr}
\item $\exists !z\in W\:\ell_n=z$;
%\item $h_n(x) R_n h_n(y)$ or $h_n(x) R_\gw h_n(y)$ or $h_n(x)=h_n(y)$ whenever $x<y$ and $x$ is sufficiently large;
\item $\ell_n R_{n+1} \ell_{n+1}$ or $\ell_n R_\gw \ell_{n+1}$ or $\ell_n=\ell_{n+1}$;
\item If $m<n$ then $\ell_m =\ell_n$ or $\ell_m R_\ga\ell_n$, for some $\ga\in (m,n]\cup\{\gw\}$.
\end{enumr}
\el

\bp\ Statement (i) is proved by (external) induction on $n$. First,
we observe that, by the polytransitivity, the relation $R_n\cup
R_\gw$ is transitive and irreflexive, for each $n$. Now it is easy
to see that the limit of $h_0$ exists, as $h_0$ only moves along
$R_0\cup R_\gw$. Suppose $\ell_{n-1}$ exists. As soon as $h_{n-1}$
reaches its limit, $h_n$ can only move along $R_n\cup R_\gw$. Hence,
$\ell_n$ exists.

Statement (ii) follows from the same consideration and the fact that
$h_{n+1}$ has to visit $\ell_n$ on its way to the limit. Statement
(iii) is obtained from (ii) by induction on $n$. \ep

\bl For all $n$, $\nat\models (\ell_n=0)$. \el

\bp\ By Lemma \ref{incr}, for all $n>M$, either $\ell_n R_\gw
\ell_{n+1}$ or $\ell_{n+1}=\ell_n$, as $R_k=\emptyset$ for $k>M$.
Since $R_\gw$ is transitive and irreflexive, there is a $z\in W$ and
an $m$ such that $\ell_n=z$, for all $n\geq m$. Assume $z\neq 0$ and
let $m$ be the minimal $n$ such that $\ell_n=z$. Then the function
$h_m$ has to come to $z$ by the second clause of
Definition~\ref{Sol}. Hence, for some $n\geq \max(M,m)$,
$\Box_m(\ell_n\neq \bar z)$. Since $\PA$ is sound, $\ell_n\neq \bar
z$ is true, which is not the case since $n\geq m$. \ep

For any modal formula $C$, let $\ell_n\fc C$ denote
$\bigvee\{\ell_n=\bar a: a\fc C\}$.

\bl \label{down} For any formula $C$, for all $n>m\geq M$,
$$\PA\vdash(\ell_n\fc C) \to (\ell_m\fc C).$$ \el

\bp\ Each $h_n$ for $n>M$ can only follow the $R_\gw$ relation. By
the persistence of $\cW$, the truth of any formula is inherited
downwards along $R_\gw$. Hence, the claim follows from Lemma
\ref{incr}. \ep Lemmas \ref{incr}, \ref{down} are obviously
formalizable in $\PA$ (uniformly in $m,n$).

Suppose $\{\phi_i:i\in I\}$ is a primitive recursive set of
formulas. With a primitive recursive program computing this set we
associate a numeration for the theory $\PA+\{\phi_i:i\in I\}$ that
will be denoted $[\phi_i: i\in I]$. We write $\Con_n[\phi_i: i\in
I]$ for $\Con_n([\phi_i: i\in I])$. In particular, if this set is a
singleton $\{\phi\}$, the formula $\Con_n[\phi]$ means the same as
$\Con_n(\ul\phi)$ and is $\PA$-equivalent to $\neg\Box_n\neg \phi$.

Using this notation, we interpret each propositional variable $p$ as
follows: $$p^*:= [\ell_n\fc p: n \geq M].$$

We prove the following two main lemmas.

\bl \label{main1} For any formula $C\in\Phi$,
$$[\ell_n\fc C: n \geq M]\vdash_{\PA} C^*.\eqno (*)$$
\el

\bp\ Induction on the build-up of $C$. The cases of propositional
variables, $\top$ and $\land$ are easy.

Assume $C=mD$ for $m<\gw$. Since $(mD)^*=\ul{\Con}_m(D^*)$ numerates
a finite extension of $\PA$, it will be sufficient in this case to
infer $\Con_m(D^*)$ from $\ell_M\fc mD$ in $\PA$. We have, by the
IH,
$$\PA+\Con_m[\ell_n\fc D:n\geq M]\vdash\Con_m(D^*).$$
However, $\Con_m[\ell_n\fc D:n\geq M]$ is equivalent to the formula
$$\textstyle \al{n\geq M} \Con_m[\bigwedge_{k=M}^n (\ell_k\fc D)],$$
which is by Lemma \ref{down} equivalent to \beq \label{conM}
\textstyle \al{n\geq M} \Con_m[\ell_n\fc D].\eeq

Thus, we are going to infer sentence \refeq{conM} from $\ell_M\fc
mD$ by formalizing the following argument in $\PA$. Assume
$\ell_M\fc mD$, hence there is a $z$ such that $z\fc D$ and $\ell_M
R_m z$. Consider the point $\ell_m$. By Lemma \ref{incr}, either
$\ell_m R_k \ell_{M}$ for some $k>m$, or $\ell_m R_\gw \ell_M$ or
$\ell_m=\ell_M$.  In each case, $\ell_m R_m z$, as $\cW$ is an
RJ-frame.

Assume $\ex{n\geq M}\neg \Con_m[\ell_n\fc D]$, then
$\Box_m(\ell_n\nfc D)$. Since (provably) $z\fc D$, we have
$\Pr_m[\ell_n\neq \bar{z}]$. Let $x_0$ be such that $\al{x\geq x_0}
h_m(x)=\ell_m$. There is a $y>x$ such that $\Prf_m(\gn{\ell_n\neq
\bar{z}},y)$. Then, $h_m(y+1)$ has to be different from $\ell_m$, a
contradiction. This shows that $\ell_M\fc mD$ implies \refeq{conM},
as required.

\bs Consider the case $C=\gw D$. Firstly, we have
$$[\Con_n(D^*):n\geq M]\vdash_\PA (\gw D)^*,$$ where we may
restrict the left hand side to $n\geq M$, since the strength of the
formulas $\Con_n$ increases with $n$. By Lemma \ref{down} and the
IH, as before,
\begin{eqnarray*} [\al{k\geq \bar n}\Con_n[\ell_k\fc
D]: n\geq M] & \vdash_\PA & [\al{k\geq
 M}\Con_n[\ell_k\fc D]: n\geq M] \\
& \vdash_\PA & [\Con_n(D^*):n\geq M].
\end{eqnarray*}
We are going to show that
$$[\ell_n\fc \gw D: n\geq M]\vdash_\PA
[\al{k\geq \bar n}\Con_n[\ell_k\fc D]: n\geq M].$$ To this end it is
sufficient to prove by an argument formalizable in $\PA$ that, for
any $n\geq M$,
$$\PA + (\ell_n\fc \gw D) \vdash \al{k\geq \bar n} \Con_n[\ell_k\fc D].$$

Consider any $n\geq M$ and assume $\ell_n\fc \gw D$. There is a
$z\in \cW$ such that $z\fc D$ and $\ell_n R_\gw z$. If $\ex{k\geq
\bar n}\Box_n(\ell_k\nfc D)$ then $\ex{k\geq \bar n}
\Box_n(\ell_k\neq \bar z)$. This means that $h_n$ must take on a
value other than $\ell_n$, a contradiction. \ep

\bl\label{main2} For any formula $C\in\Phi$,
$$\ul{\ell_0\neq 0}\lor C^*\vdash_{\PA} [\ell_n\fc C: n \geq M].\eqno (**)$$
\el

\bp\  Induction on the build-up of $C$. The cases of propositional
variables, $\top$ and $\land$ are trivial.

Assume $C=mD$ for $m<\gw$. Since $\ell_0\neq 0$ is equivalent to a
$\Sigma_1$-formula, we obtain
\begin{eqnarray*} \PA + \ell_0\neq
0\land \Con_m(D^*) & \vdash & \Box_0(\ell_0\neq 0) \land \Con_m(D^*)
\\ & \vdash & \Con_m(\ul{\ell_0\neq 0} \lor D^*)
\\ & \vdash & \Con_m[\ell_k\fc D: k \geq M]
\\ & \vdash & \al{k\geq M}\Con_m[\ell_k\fc D].
\end{eqnarray*}

Thus, it is sufficient to prove, for each $n\geq M$, that
$$\PA\vdash \ell_0\neq 0\land \ell_n\nfc mD \to \ex{k\geq M}
\Box_m(\ell_k\nfc D).$$ We reason as follows.

Assume $\ell_n\nfc mD$. By Lemma \ref{incr} we have $\ell_m  R_{k}
\ell_n$, for some $k>m$, or $\ell_m R_\gw \ell_n$ or
$\ell_m=\ell_n$. Since $\cW$ is a RJ-frame, in each case
\beq\ell_m\nfc mD. \label{t} \eeq
%This means that $\al{z}(a R_m z \to z\nfc D)$.
Let $a:=\ell_m$, we have $$\ex{x}(h_m(x)=a \land \al{y\geq
x}(h_{m-1}(y)=h_{m-1}(x))).\eqno (L_m(a))$$ (After $h_m$ attains its
limit $a$, the function $h_{m-1}$ only has the possibility to make a
single move to $a$, if it moves at all.) The statement $L_m(a)$ is
expressible by a $\Sigma_{m+1}$-formula. Hence, $\Box_m L_m(\bar
a)$. Moreover, $L_m(a)$ implies that $h_m$ goes along the $R_m\cup
R_\gw$ relations from $a$ onwards. Hence, for any $k\geq m$,
$L_m(a)$ implies $\ell_k\in R^*_m(a)\cup\{a\}$, where
$$R^*_m(a):=\{y\in W: \ex{\ga\geq m} a R_\ga y\}.$$
Therefore, we conclude: \beq \label{t2} \al{k\geq m}
\Box_m(\ell_k\in R^*_m(\bar a)\cup\{\bar a\}). \eeq

On the other hand, we claim that $\al{z\in R^*_m(a)} z\nfc D$.
Indeed, if $a R_\ga z$, $\ga\geq m$ and $z\fc D$, then $a\fc \ga D$
whence $a\fc mD$, by the $\Phi$-monotonicity of $\cW$. This
contradicts \refeq{t}. Since the formula $\al{z\in R^*_m(a)} z\nfc
D$ is bounded, it follows that \beq \label{t3} \Box_m(\al{z\in
R^*_m(\bar a)} z\nfc D).\eeq

Next we consider the minimal $i\leq m$ such that $\ell_i=\ell_m=a$.
Since we assume $\ell_0\neq 0$, we also have $\ell_m=a\neq 0$. The
function $h_i$ could only have come to $a$ by the second clause of
the definition of $h_i$, therefore we obtain
$$\ex{k\geq\max(i,M)}\Box_i(\ell_k\neq \bar a).$$
Since $i\leq m$, obviously we can weaken this to \beq \ex{k\geq M}
\Box_m(\ell_k\neq \bar a). \eeq Together with \refeq{t2} this yields
$\ex{k\geq M} \Box_m(\ell_k\in R^*_m(\bar a)),$ whence by \refeq{t3}
we obtain $\ex{k\geq M} \Box_m(\ell_k\nfc D),$ as required.

\bigskip
Finally, consider the case $C=\gw D$. We have $$(\gw D)^*
=[\Con_n(D^*):n\in\gw].$$ For each of the axioms $\Con_n(D^*)$ such
that $n\geq M$, by the IH and persistence we obtain, as before,
$$\PA+\ell_0\neq 0\land \Con_n(D^*) \vdash \al{k\geq M}\Con_n[\ell_k\fc D].$$
This fact is formalizable in $\PA$ uniformly in $n$, therefore
$$\ul{\ell_0\neq 0}\lor (\gw D)^*\vdash_{\PA}
[\al{k\geq M}\Con_n[\ell_k\fc D]:n\geq M].$$ We are going to show
that
$$\ul{\ell_0\neq 0}\lor [\al{k\geq M}\Con_n[\ell_k\fc D]:n\geq M]
\vdash_\PA [\ell_n\fc \gw D: n>M],$$ which completes the proof,
because by Lemma \ref{down} $$[\ell_n\fc \gw D: n> M]\vdash_\PA
[\ell_n\fc \gw D: n\geq M].$$

Thus, we prove by an argument formalizable in $\PA$ uniformly in $n$
that, for any $n>M$,
$$\PA +(\ell_0\neq 0\land\ell_n\nfc \gw D)\vdash
\ex{k\geq M} \Box_n(\ell_k\nfc D).$$ Assume $n>M$, $\ell_n\nfc \gw
D$ and let $a:=\ell_n$. Since $\ell_0\neq 0$ we have $a\neq 0$.
Consider the minimal $m\leq n$ such that $\ell_m=\ell_n=a$. As
before, we first show that \beq \label{a} \al{k\geq
\max(m,M)}\Box_n(\ell_k\in R_\gw(\bar a)\cup\{\bar a\}).\eeq

We consider two cases. If $m>M$ then we have both
 $L_m(a)$ and $\Box_m(L_m(\bar a))$.
Moreover, for $k\geq m>M$, from $L_m(\bar{a})$ we can infer
$\ell_k\in R_\gw(\bar a)\cup\{\bar a\}$, since $h_k$ can only make
moves along the $R_\gw$ relation from $a$ onwards. (Here,
$R_\gw(a):=\{x\in \cW:a R_\gw x\}$.) Hence, $\al{k\geq
m}\Box_m(\ell_k\in R_\gw(\bar a)\cup\{\bar a\})$ and the claim
holds.

If $m\leq M<n$ we first notice that $\ell_M=a$. Then,
$\Box_n(\ell_M=\bar a)$, since the formula $\ell_M=\bar a$ is
$\Sigma_{M+2}\subseteq \Sigma_{n+1}$. Moreover, $\ell_M=a$ implies
$\ell_k\in R_\gw(a)\cup\{a\}$, for any $k>M$, since $h_k$ will only
be able to move along $R_\gw$ from $a$ onwards. Thus, we obtain
$\al{k\geq M}\Box_n(\ell_k\in R_\gw(\bar a)\cup\{\bar a\}),$ and the
claim also holds.

Secondly, we note that $\ell_k\in R_\gw(a)$ implies $\ell_k\nfc D$,
as $\al{z\in R_\gw(a)} z\nfc D$. Hence, by \refeq{a},
\beq\label{b}\al{k\geq \max(m,M)}\Box_n(\ell_k\nfc D \lor
\ell_k=\bar a).\eeq

Thirdly, since $a\neq 0$, by the definition of $h_m$ we have
$$\ex{k\geq \max(M,m)}\Box_m(\ell_k\neq \bar a).$$ Together with \refeq{b}
this yields $\ex{k\geq \max(M,m)}\Box_n(\ell_k\nfc D)$. In either
case we obtain $\ex{k\geq M}\Box_n(\ell_k\nfc D)$. \ep

Recall that at the node $1\in \cW$ there holds $1\fc A$ and $1\nfc
B$. Let $\gs$ denote $ [\ell_n=1 :n\geq M]$ and $S$ denote the
corresponding theory. By Lemma \ref{main1}, $\gs\vdash_\PA A^*.$ On
the other hand, by Lemma \ref{main2}
$$\ul{\ell_0\neq 0}\lor B^*\vdash_\PA [\ell_n\fc B:n\geq M]
\vdash_\PA [\ell_n\neq 1:n\geq M].$$ Hence, $A^*\vdash B^*$ yields
$S\vdash \ell_M\neq 1.$ It follows that $S$ is inconsistent. Since
$\PA\vdash \ell_n=1\to \ell_m=1,$ for all $m\leq n$, the
inconsistency of $S$ yields a $\PA$-proof of $\ell_n\neq 1$, for
some $n\geq M$. This means that $h_0$ must eventually take on a
value other than $0$, hence $\ell_0\neq 0$. But this is impossible,
since $\ell_0=0$ is true in the standard model. \ep

\bex By Example \ref{ex},  $\gw p\land\gw q\nvdash_\RCo \gw(p\land
q)$. By Theorem \ref{arcomp}, this means that there are theories
$S,T$ containing $\PA$ such that
$$\PA+\RFN(S)+\RFN(T)\nvdash
\RFN(S+T).$$ Primitive recursive numerations of these theories can
be obtained from the proof of Theorem \ref{arcomp} applied to (an
irreflexive version of) the three-element Kripke model described in
Example \ref{ex}.

We remark that none of these two theories can have bounded
arithmetical complexity over $\PA$. Suppose $S$ is axiomatized by a
set of $\Pi_{n+1}$-sentences over $\PA$. Then, by Lemma \ref{reflc},
$\ul{\Con}_n(\gs)\vdash_\PA \gs$. By Corollary \ref{j-ax}, it
follows that
$$\PA+\Con_\gw(\gs)+\Con_\gw(\tau)\vdash \Con_\gw(\tau\lor
\ul{\Con}_n(\gs))\vdash \Con_\gw(\tau\lor \gs).$$ This shows that
the use of infinitely axiomatized theories to interpret
propositional variables is necessary for the validity of Theorem
\ref{arcomp}.

\eex

\section{Conclusions}

We believe that positive provability logic, despite the absence of
L\"ob's axiom, strikes a good balance between expressivity and
efficiency (the latter can be understood formally, as the
computational efficiency, as well as informally, in the sense of
convenience). Together with \cite{Das12} this paper shows that
positive logic can be nicely treated both syntactically and
semantically. More importantly, it has very natural proof-theoretic
interpretations not extendable to the full modal logic language.

There are many questions related to this logic that can be further
investigated. One direction is to study normal positive logics along
the lines of the usual normal modal logics. In particular, we are
interested in their efficient proof systems, general results on
axiomatization and completeness, interpolation properties, and so
on.

Another direction is the study of different arithmetical
interpretations of positive provability logic. For example, one can
consider from this point of view transfinite iterations of
consistency assertions (or of higher reflection principles). That
is, one can introduce modalities $\Diamond^\ga$, for each ordinal
$\ga$ of some canonical ordinal notation system, and interpret them
as the schemata $\Con^\ga$ related to  the so-called Turing
progressions: $\Con^0(\gs)=\Con(\gs);$
$\Con^\ga(\gs)=\{\Con[\Con^\gb(\gs)]:\gb<\ga\}$. It would be
interesting to find a complete axiomatization of the corresponding
positive logic.

Another generalization is to consider stronger reflection schemata
definable in the extensions of arithmetical language, e.g., in the
second order arithmetic or in the arithmetic enriched by
truthpredicates. This generalization is particularly interesting
from the point of view of applications in the ordinal analysis of
predicative theories.

\section{Acknowledgements}

This paper is dedicated to Sergei N. Artemov who guided the author
into the fascinating world of provability logic some 25 years ago.

Thanks are due to E. Dashkov, V. Krupski, and I. Shapirovsky for
their helpful comments and suggestions for improvement. This work
was supported by the Russian Foundation for Basic Research, Russian
Presidential Council for Support of Leading Scientific Schools, and
the STCP-CH-RU project ``Computational proof theory''.

\bibliographystyle{model1a-num-names}
%\bibliography{ref-all2}

\appendix
\section{Arithmetization of the Solovay functions}

Given a finite irreflexive RJ-frame $\cW$ we would like to build a
family of arithmetical functions $h_n:\gw\to \cW$, for all
$n\in\gw$, satisfying the following conditions: \bi
\item The graph of each $h_n$ is definable by a
formula $H_n$ which is $\Delta_{n+1}$ in $\PA$;
\item There is a primitive recursive function $\phi_e:n\mapsto
\gn{H_n}$, where $e$ denotes the primitive recursive index of this
function;
\item Each $h_n$ provably satisfies the clauses of Definition
\ref{Sol}. \ei

These objects will be constructed using the formalized recursion
theorem. The main unknown is the index $e$.

Firstly, we stipulate that the limit statements $\ell_n=z$ are
abbreviations for the formulas $\ex{N}\al{k>N} H_n(k,z)$. Secondly,
we fix a primitive recursive function $g_0$ such that
$$\gn{\ell_n=\bar z} = g_0(\gn{H_n},z) = g_0(\phi_e(n),z).$$
We see that the function $g(e,n,z):=g_0(\phi_e(n),z)$ is provably
total recursive in $\PA$, hence it is definable by an arithmetical
$\Delta_1$-formula.

Using $g$, Definition \ref{Sol} can be rewritten to define the
graphs of $h_n$ in the language of arithmetic with the unknown $e$
as an extra parameter. We denote such parametrized versions of the
formulas $H_n$ by $H'_n$. Each formula $H'_n$ uses the formulas
$H'_0,\dots,H'_{n-1}$ as subformulas to express the first clause of
the definition of $h_n$. Thus, we obtain a sequence of formulas of
the following form:

\begin{eqnarray*}
H'_0(e,x,y) & \eqv & A_0(e,x,y) \\
H'_1(e,x,y) & \eqv & A_1(H'_0;e,x,y) \\
&\dots & \\
H'_n(e,x,y) & \eqv & A_n(H'_0,\dots,H'_{n-1};e,x,y)
\end{eqnarray*}
Here, the formulas $A_n$ directly mimic Definition \ref{Sol}. It is
easy to convince oneself that the arithmetical complexity of the
formulas $A_n$ (and hence, of the formulas $H'_n$) is $\Delta_{n+1}$
in $\PA$. The most complex part of the definition is the formula
$\ex{k\geq \max(M,n)}\Prf_n(\gn{\ell_k\neq \bar{z}},x)$ occurring in
the second clause. In the formula $A_n$ this part takes the form
$$\ex{k\geq \max(M,n)}\Prf_n(g(e,k,z),x).$$
Observe that the predicate $\Prf_n$ is $\Delta_{n+1}$ (even
$\Delta_0(\Sigma_n)$), and the existential quantifier $\ex{k}$ can,
in fact, be bounded by $x$. This yields a $\Delta_{n+1}$-formula.

It is also clear that each formula $H'_n$ is obtained from the
previous ones in a primitive recursive way. Therefore, there is a
primitive recursive function $F$ satisfying $F(e,n)= \gn{H'_n(\bar
e,x,y)}$. Finally, we obtain the required number $e$ by applying (a
formalized version of) recursion theorem for primitive recursive
functions: $\phi_e(n):= F(e,n).$ Then we can define $H_n(x,y):=
H_n'(\bar e,x,y)$.

\end{document}